\documentclass{amsart}
\usepackage{amscd, amssymb, amsfonts}
\usepackage{amsmath,amsthm}
\usepackage{amssymb}

  [2002/01/19 v2.2g %
    AMS font definitions%
  ]
\DeclareFontFamily{U}{euf}{}
\DeclareFontShape{U}{euf}{m}{n}{%
  <5><6><7><8><9>gen*eufm%
  <10><10.95><12><14.4><17.28><20.74><24.88>eufm10%
  }{}
\DeclareFontShape{U}{euf}{b}{n}{%
  <5><6><7><8><9>gen*eufb%
  <10><10.95><12><14.4><17.28><20.74><24.88>eufb10%
  }{}

  [2002/01/19 v2.2g %
    AMS font definitions%
  ]
\DeclareFontFamily{U}{msb}{}
\DeclareFontShape{U}{msb}{m}{n}{%
  <5><6><7><8><9>gen*msbm%
  <10><10.95><12><14.4><17.28><20.74><24.88>msbm10%
  }{}

  [2002/01/19 v2.2g %
    AMS font definitions%
  ]
\DeclareFontFamily{U}{msa}{}
\DeclareFontShape{U}{msa}{m}{n}{%
  <5><6><7><8><9>gen*msam%
  <10><10.95><12><14.4><17.28><20.74><24.88>msam10%
  }{}

\newtheorem{theorem}{Theorem}[section]
\newtheorem{lemma}[theorem]{Lemma}

\newtheorem{corollary}[theorem]{Corollary}

\theoremstyle{definition}

\theoremstyle{remark}

\numberwithin{equation}{section}

\begin{document}

\title[On $q$-Euler numbers related to the modified $q$-Bernstein polynomials]
{On $q$-Euler numbers related to the modified $q$-Bernstein polynomials}

\author{Min-Soo Kim, Daeyeoul Kim and Taekyun Kim}

\begin{abstract}
We consider $q$-Euler numbers and polynomials and $q$-Stirling numbers of first and second kinds.
Finally, we investigate some interesting properties of the modified $q$-Bernstein polynomials related to $q$-Euler numbers and
$q$-Stirling numbers
by using fermionic $p$-adic integrals on $\mathbb Z_p.$
\end{abstract}

\address{Department of Mathematics, KAIST, 373-1 Guseong-dong, Yuseong-gu, Daejeon 305-701, S. Korea}
\email{minsookim@kaist.ac.kr}

\address{National Institute for Mathematical Sciences, Doryong-dong, Yuseong-gu, Daejeon 305-340, Republic of KOREA}
\email{daeyeoul@nims.re.kr}

\address{Division of General Education-Mathematics, Kwangwoon University, Seoul, 139-701, S. Korea}
\email{tkkim@kw.ac.kr}



\maketitle



\def\ord{\text{ord}_p}
\def\C{\mathbb C_p}
\def\BZ{\mathbb Z}
\def\Z{\mathbb Z_p}
\def\Q{\mathbb Q_p}
\def\wh{\widehat}

\section{Introduction}
\label{Intro}

Let $C[0,1]$ be the set of continuous functions on $[0,1].$ The classical Bernstein polynomials of degree $n$
for $f\in C[0,1]$ are defined by
\begin{equation}\label{cla-def-ber}
\mathbb B_{n}(f)=\sum_{k=0}^nf\left(\frac kn\right)B_{k,n}(x),\quad 0\leq x\leq1
\end{equation}
where $\mathbb B_{n}(f)$ is called the Bernstein operator and
\begin{equation}\label{or-ber-poly}
B_{k,n}(x)=\binom nk x^k(x-1)^{n-k}
\end{equation}
are called the Bernstein basis polynomials (or the Bernstein polynomials of
degree $n$) (see \cite{SA}).
Recently, Acikgoz and Araci have studied the generating function for Bernstein
polynomials (see \cite{AA,AA2}). Their generating function for $B_{k,n}(x)$ is given by
\begin{equation}\label{AA-gen}
F^{(k)}(t,x)=\frac{t^ke^{(1-x)t}x^k}{k!}=\sum_{n=0}^\infty B_{k,n}(x)\frac{t^n}{n!},
\end{equation}
where $k=0,1,\ldots$ and $x\in[0,1].$ Note that
$$B_{k,n}(x)=\begin{cases}
\binom nk x^k(1-x)^{n-k}&\text{if } n\geq k \\
0,&\text{if }n<k
\end{cases}$$
for $n=0,1,\ldots$ (see \cite{AA,AA2}).

Let $p$ be an odd prime number.
Throughout this paper, $\mathbb Z_p, \mathbb Q_p$ and $\mathbb C_p$
will denote the ring of $p$-adic rational
integers, the field of $p$-adic rational numbers and the completion
of the algebraic closure of $\mathbb Q_p,$ respectively.
Let $v_p$ be the normalized exponential valuation of $\mathbb C_p$ with
$|p|_p=p^{-1}.$

Throughout this paper, we use the following notation
$$[x]_q=\frac{1-q^x}{1-q}\quad\text{and}\quad [x]_{-q}=\frac{1-(-q)^x}{1+q}$$
(cf. \cite{KT1,KT2,KT3,KT7}).
Let $\mathbb N$ be the natural numbers and $\mathbb Z_+=\mathbb N\cup\{0\}.$
Let $UD(\mathbb Z_p)$ be the space of uniformly differentiable function on $\mathbb Z_p.$

Let $q\in\C$ with $|1-q|_p<p^{-1/(p-1)}$ and $x\in\Z.$ Then $q$-Bernstein type operator for $f\in UD(\Z)$ is defined by
(see \cite{KJY,KCK})
\begin{equation}\label{B-op}
\begin{aligned}
\mathbb B_{n,q}(f)&=\sum_{k=0}^nf\left(\frac kn\right)\binom nk[x]_q^k[1-x]_q^{n-k}\\
&=\sum_{k=0}^nf\left(\frac kn\right)B_{k,n}(x,q),
\end{aligned}
\end{equation}
for $k,n\in \mathbb Z_+,$ where
\begin{equation}\label{Ber-def}
B_{k,n}(x,q)=\binom nk[x]_q^k[1-x]_q^{n-k}
\end{equation}
is called the modified $q$-Bernstein polynomials
of degree $n.$
When we put $q\to1$ in (\ref{Ber-def}), $[x]_q^k\to x^k,[1-x]_q^{n-k}\to(1-x)^{n-k}$
and we obtain the classical Bernstein polynomial, defined by (\ref{or-ber-poly}).
We can deduce very easily from (\ref{Ber-def}) that
\begin{equation}\label{recu-ber}
B_{k,n}(x,q)=[1-x]_qB_{k,n-1}(x,q)+[x]_qB_{k-1,n-1}(x,q)
\end{equation}
(see \cite{KJY}).
For $0 \leq k \leq n,$ derivatives of the $n$th degree modified $q$-Bernstein polynomials are polynomials of
degree $n-1:$
\begin{equation}\label{der-Ber-def}
\frac{\text{d}}{\text{d}x}B_{k,n}(x,q)=n(q^xB_{k-1,n-1}(x,q)-q^{1-x}B_{k,n-1}(x,q))\frac{\ln q}{q-1}
\end{equation}
(see \cite{KJY}).

The Bernstein polynomials can also be defined in many
different ways. Thus, recently, many applications of these polynomials have been looked
for by many authors.
In recent years, the $q$-Bernstein polynomials have been investigated and studied by many authors in many
different ways (see \cite{KJY,KCK,SA} and references therein \cite{GG,Ph}).
In \cite{Ph}, Phillips gave many results concerning the $q$-integers, and an account
of the properties of $q$-Bernstein polynomials. He gave many applications of these polynomials
on approximation theory.
In \cite{AA,AA2}, Acikgoz and Araci have introduced several type Bernstein polynomials. The Acikgoz and Araci paper to announce in the
conference is actually motivated to write this paper.
In \cite{SA}, Simsek and Acikgoz constructed a new generating function of the $q$-Bernstein type polynomials
and established elementary properties of this function.
In \cite{KJY}, Kim, Jang and Yi proposed the modified $q$-Bernstein polynomials of degree $n,$ which
are different $q$-Bernstein polynomials of Phillips.
In \cite{KCK}, Kim, Choi and Kim investigated some interesting properties of the modified $q$-Bernstein polynomials of degree $n$
related to $q$-Stirling numbers and Carlitz's $q$-Bernoulli numbers.

In the present paper,
we consider $q$-Euler numbers, polynomials and $q$-Stirling numbers of first and second kinds.
We also investigate some interesting properties of the modified $q$-Bernstein polynomials of degree $n$
related to $q$-Euler numbers and
$q$-Stirling numbers by using fermionic $p$-adic integrals on $\mathbb Z_p.$

\section{$q$-Euler numbers and polynomials
related to the fermionic $p$-adic integrals on $\mathbb Z_p$}

For $N\geq1,$ the fermionic $q$-extension $\mu_q$
of the $p$-adic Haar distribution $\mu_{\text{Haar}}:$
\begin{equation}\label{mu}
\mu_{-q}(a+p^N\Z)=\frac{(-q)^a}{[p^{N}]_{-q}}
\end{equation}
is known as a measure on $\Z,$ where $a+p^N\Z=\{ x\in\Q\mid |x-a|_p\leq p^{-N}\}$ (cf. \cite{KT1,KT4}).
We shall write $d\mu_{-q}(x)$ to remind ourselves that $x$ is the variable
of integration.
Let $UD(\mathbb Z_p)$ be the space of uniformly differentiable function on $\mathbb Z_p.$
Then $\mu_{-q}$ yields the fermionic $p$-adic $q$-integral of a function $f\in UD(\mathbb Z_p):$
\begin{equation}\label{Iqf}
I_{-q}(f)=\int_{\Z} f(x)d\mu_{-q}(x)=\lim_{N\rightarrow\infty}\frac{1+q}{1+q^{p^N}}
\sum_{x=0}^{p^N-1}f(x)(-q)^x
\end{equation}
(cf. \cite{KT4,KT5,KT6,KT8}).
Many interesting properties of (\ref{Iqf}) were studied by many authors
(see \cite{KT4,KT5} and the references given there).
Using (\ref{Iqf}), we have the fermionic $p$-adic invariant integral on $\mathbb Z_p$ as
follows:
\begin{equation}\label{de-2}
\lim_{q\to-1}I_q(f)=I_{-1}(f)=\int_{\mathbb Z_p}f(a)d\mu_{-1}(x).
\end{equation}
For $n\in\mathbb N,$ write $f_n(x)=f(x+n).$ We have
\begin{equation}\label{de-3}
I_{-1}(f_n)=(-1)^nI_{-1}(f)+2\sum_{l=0}^{n-1}(-1)^{n-l-1}f(l).
\end{equation}
This identity is obtained by Kim in \cite{KT4}
to derives interesting properties and relationships involving $q$-Euler numbers and polynomials.
For $n\in\mathbb Z_+,$ we note that
\begin{equation}\label{q-Euler-numb}
I_{-1}([x]_q^n)=\int_{\Z} [x]_q^n d\mu_{-1}(x)=E_{n,q},
\end{equation}
where $E_{n,q}$ are the $q$-Euler numbers (see \cite{KT09}).
It is easy to see that $E_{0,q}=1.$ For $n\in\mathbb N,$ we have
\begin{equation}\label{q-Eu-rec}
\begin{aligned}
\sum_{l=0}^n\binom nlq^lE_{l,q}&=\sum_{l=0}^n\binom nlq^l \lim_{N\rightarrow\infty}\sum_{x=0}^{p^N-1}[x]_q^l(-1)^x \\
&=\lim_{N\rightarrow\infty}\sum_{x=0}^{p^N-1}(-1)^x(q[x]_q+1)^n \\
&=\lim_{N\rightarrow\infty}\sum_{x=0}^{p^N-1}(-1)^x[x+1]_q^{n} \\
&=-\lim_{N\rightarrow\infty}\sum_{x=0}^{p^N-1}(-1)^x([x]_q^{n}+[p^N]_q^n) \\
&=-E_{n,q}.
\end{aligned}
\end{equation}
From this formula, we have the following recurrence formula
\begin{equation}\label{q-Euler-recur}
E_{0,q}=1,\qquad (qE+1)^n+E_{n,q}=0\quad\text{if }n\in\mathbb N
\end{equation}
with the usual convention of replacing $E^l$ by $E_{l,q}.$
By the simple calculation of the fermionic $p$-adic invariant integral on $\mathbb Z_p,$
we see that
\begin{equation}\label{q-E-ex}
E_{n,q}=\frac{2}{(1-q)^n}\sum_{l=0}^n\binom nl(-1)^l\frac{1}{1+q^l},
\end{equation}
where $\binom nl=n!/l!(n-l)!=n(n-1)\cdots(n-l+1)/l!.$
Now, by introducing the following equations:
\begin{equation}\label{intr-q-ex}
[x]_{\frac1q}^n=q^{n}q^{-nx}[x]_q^n\quad\text{and}\quad q^{-nx}=\sum_{m=0}^\infty(1-q)^m\binom{n+m-1}{m}[x]_q^m
\end{equation}
into (\ref{q-Euler-numb}), we find that
\begin{equation}\label{q-Eu-inv}
E_{n,\frac1q}=q^{n}\sum_{m=0}^\infty(1-q)^m\binom{n+m-1}{m}E_{n+m,q}.
\end{equation}
This identity is a peculiarity of the $p$-adic $q$-Euler numbers, and the classical Euler numbers do not seem
to have a similar relation.
Let $F_q(t)$ be the generating function of the $q$-Euler numbers. Then we obtain
\begin{equation}\label{e-E-gen}
\begin{aligned}
F_q(t)&=\sum_{n=0}^\infty E_{n,q}\frac{t^n}{n!} \\
&=\sum_{n=0}^\infty\frac{2}{(1-q)^n}\sum_{l=0}^n(-1)^l\binom nl\frac{1}{1+q^l}\frac{t^n}{n!}\\
&=2e^{\frac{t}{1-q}}\sum_{k=0}^\infty\frac{(-1)^k}{(1-q)^k}\frac{1}{1+q^k}\frac{t^k}{k!}.
\end{aligned}
\end{equation}
From (\ref{e-E-gen}) we note that
\begin{equation}\label{e-E-gen-mod}
F_q(t)=2e^{\frac{t}{1-q}}\sum_{n=0}^\infty(-1)^ne^{\left(\frac{-q^n}{1-q}\right)t}
=2\sum_{n=0}^\infty(-1)^ne^{[n]_qt}.
\end{equation}
It is well-known that
\begin{equation}\label{q-Euler}
I_{-1}([x+y]^n)=\int_{\Z} [x+y]^n d\mu_{-1}(y)=E_{n,q}(x),
\end{equation}
where $E_{n,q}(x)$ are the $q$-Euler polynomials (see \cite{KT09}).
In the special case $x=0,$ the numbers $E_{n,q}(0)=E_{n,q}$ are referred to as the $q$-Euler numbers.
Thus we have
\begin{equation}\label{E-p-n}
\begin{aligned}
\int_{\Z} [x+y]^n d\mu_{-1}(y)&=\sum_{k=0}^n\binom nk[x]_q^{n-k}q^{kx}\int_{\Z} [y]^k d\mu_{-1}(y) \\
&=\sum_{k=0}^n\binom nk[x]_q^{n-k}q^{kx}E_{k,q} \\
&=(q^xE+[x]_q)^n.
\end{aligned}
\end{equation}
It is easily verified, using (\ref{e-E-gen-mod}) and (\ref{q-Euler}), that
the $q$-Euler polynomials $E_{n,q}(x)$ satisfy the following formula:
\begin{equation}\label{e-p-exp}
\begin{aligned}
\sum_{n=0}^\infty E_{n,q}(x)\frac{t^n}{n!}&=\int_{\Z}e^{[x+y]_qt}d\mu_{-1}(y) \\
&=\sum_{n=0}^\infty\frac{2}{(1-q)^n}\sum_{l=0}^n(-1)^l\binom nl\frac{q^{lx}}{1+q^l}\frac{t^n}{n!}\\
&=2\sum_{n=0}^\infty(-1)^ne^{[n+x]_qt}.
\end{aligned}
\end{equation}
Using formula (\ref{e-p-exp}) when $q$ tends to 1, we can readily derive the Euler polynomials, $E_n(x),$ namely,
$$\int_{\Z}e^{(x+y)t}d\mu_{-1}(y)=\frac{2e^{xt}}{e^t+1}=\sum_{n=0}^\infty E_n(x)\frac{t^n}{n!}$$
(see \cite{KT4}).
Note that $E_n(0)=E_n$ are referred to as the $n$th Euler numbers.
Comparing the coefficients of ${t^n}/{n!}$ on both sides of (\ref{e-p-exp}), we have
\begin{equation}\label{q-Euler-pro}
E_{n,q}(x)=2\sum_{m=0}^\infty(-1)^m[m+x]_q^n=\frac{2}{(1-q)^n}\sum_{l=0}^n(-1)^l\binom nl\frac{q^{lx}}{1+q^l}.
\end{equation}

We refer to $[n]_q$ as a $q$-integer and note that $[n]_q$ is a continuous function of $q.$ In an obvious way we also
define a $q$-factorial,
$$[n]_q!=\begin{cases}[n]_q[n-1]_q\cdots[1]_q&n\in\mathbb N, \\
1,&n=0
\end{cases}$$
and a $q$-analogue of binomial coefficient
\begin{equation}\label{q-binom}
\binom xn_q=\frac{[x]_q!}{[x-n]_q![n]_q!}=\frac{[x]_q[x-1]_q\cdots[x-n+1]_q}{[n]_q!}
\end{equation}
 (cf. \cite{KT6,KT09}).
Note that
$$\lim_{q\to1}\binom xn_q=\binom xn=\frac{x(x-1)\cdots(x-n+1)}{n!}.$$
It readily follows from (\ref{q-binom}) that
\begin{equation}\label{q-binom-1}
\binom xn_q=\frac{(1-q)^nq^{-\binom n2}}{[n]_q!}\sum_{i=0}^nq^{\binom i2}\binom ni_q(-1)^{n+i}q^{(n-i)x}
\end{equation}
(cf. \cite{KT09,KT7}).
It can be readily seen that
\begin{equation}\label{q-id-1}
q^{lx}=([x]_q(q-1)+1)^l=\sum_{m=0}^l\binom lm(q-1)^m[x]_q^m.
\end{equation}
Thus by (\ref{q-Euler}) and (\ref{q-id-1}), we have
\begin{equation}\label{q-binom-int}
\int_{\Z}\binom xn_qd\mu_{-1}(x)=\frac{(q-1)^n}{[n]_q!q^{\binom n2}}\sum_{i=0}^nq^{\binom i2}\binom ni_q(-1)^{i}
\sum_{j=0}^{n-i}\binom{n-i}j(q-1)^jE_{j,q}.
\end{equation}
From now on, we use the following notation
\begin{equation}\label{1-st-ca}
\frac{[x]_q!}{[x-k]_{q}!}=q^{-\binom k2}\sum_{l=0}^ks_{1,q}(k,l)[x]_q^l,\quad k\in \mathbb Z_+,
\end{equation}
\begin{equation}\label{2-st-ca}
[x]_q^n=\sum_{k=0}^nq^{\binom k2}s_{2,q}(n,k)\frac{[x]_q!}{[x-k]_{q}!},\quad n\in \mathbb Z_+
\end{equation}
(see \cite{KT7}).
From (\ref{1-st-ca}), (\ref{2-st-ca}) and (\ref{q-id-1}), we calculate the following consequence
\begin{equation}\label{2-st-eq}
\begin{aligned}
\;[x]_q^n&=\sum_{k=0}^nq^{\binom k2}s_{2,q}(n,k)\frac1{(1-q)^k}\sum_{l=0}^k\binom kl_q
q^{\binom l2}(-1)^lq^{l(x-k+1)} \\
&=\sum_{k=0}^nq^{\binom k2}s_{2,q}(n,k)\frac1{(1-q)^k}\sum_{l=0}^k\binom kl_q q^{\binom l2+l(1-k)}(-1)^l \\
&\quad\times\sum_{m=0}^l\binom lm(q-1)^m[x]_q^m
 \\
 &=\sum_{k=0}^nq^{\binom k2}s_{2,q}(n,k)\frac1{(1-q)^k} \\
&\quad\times\sum_{m=0}^k(q-1)^m\left(\sum_{l=m}^k\binom kl_qq^{\binom l2+l(1-k)}\binom lm(-1)^l\right)[x]_q^m.
\end{aligned}
\end{equation}
Therefore, we obtain the following theorem.

\begin{theorem}
For $n\in \mathbb Z_+,$
$$E_{n,q}=\sum_{k=0}^n\sum_{m=0}^k\sum_{l=m}^k q^{\binom k2}s_{2,q}(n,k)
(q-1)^{m-k}\binom kl_qq^{\binom l2+l(1-k)}\binom lm(-1)^{l+k}E_{m,q}.$$
\end{theorem}

By (\ref{q-id-1}) and simple calculation, we find that
\begin{equation}\label{q-b-int}
\begin{aligned}
\sum_{m=0}^n\binom nm(q-1)^mE_{m,q}&=\int_{\Z}q^{nx}d\mu_{-1}(x)\\
&=\sum_{k=0}^n(q-1)^kq^{\binom k2}\binom nk_q\int_{\Z}\prod_{i=0}^{k-1}[x-i]_{q}d\mu_{-1}(x) \\
&=\sum_{k=0}^n(q-1)^k\binom nk_q\sum_{m=0}^ks_{1,q}(k,m)\int_{\Z}[x]_q^md\mu_{-1}(x) \\
&=\sum_{m=0}^n\left(\sum_{k=m}^n(q-1)^k\binom nk_qs_{1,q}(k,m)\right)E_{m,q}.
\end{aligned}
\end{equation}
Therefore, we deduce the following theorem.

\begin{theorem}
For $n\in \mathbb Z_+,$
$$\sum_{m=0}^n\binom nm(q-1)^mE_{m,q}=\sum_{m=0}^n\sum_{k=m}^n(q-1)^k\binom nk_qs_{1,q}(k,m)E_{m,q}.$$
\end{theorem}

\begin{corollary}\label{bi-qbi}
For $m,n\in \mathbb Z_+$ with $m\leq n,$
$$\binom nm(q-1)^m=\sum_{k=m}^n(q-1)^k\binom nk_qs_{1,q}(k,m).$$
\end{corollary}

By (\ref{q-Euler-pro}) and Corollary \ref{bi-qbi}, we obtain the following corollary.

\begin{corollary}\label{E-s1-bi}
For $n\in \mathbb Z_+,$
$$E_{n,q}(x)=\frac{2}{(1-q)^n}\sum_{l=0}^n\sum_{k=l}^n(-1)^l
(q-1)^{k-l}\binom nk_qs_{1,q}(k,l)
\frac{q^{lx}}{1+q^l}.$$
\end{corollary}

It is easy to see that
\begin{equation}\label{ea-id}
\binom nk_q=\sum_{l_0+\cdots+l_k=n-k}q^{\sum_{i=0}^kil_i}
\end{equation}
(cf. \cite{KT7}). From (\ref{ea-id}) and Corollary \ref{E-s1-bi}, we can also derive the following interesting formula for
$q$-Euler polynomials.

\begin{theorem}
For $n\in \mathbb Z_+,$
$$E_{n,q}(x)=2\sum_{l=0}^n\sum_{k=l}^n\sum_{l_0+\cdots+l_k=n-k}q^{\sum_{i=0}^kil_i}\frac1{(1-q)^{n+l-k}}
s_{1,q}(k,l)(-1)^k
\frac{q^{lx}}{1+q^l}.$$
\end{theorem}

These polynomials are related to the many branches of Mathematics, for example,
combinatorics, number theory, discrete probability distributions for finding higher-order
moments (cf. \cite{KT6,KT09,KT8}). By substituting $x=0$ into
the above, we have
$$E_{n,q}=2\sum_{l=0}^n\sum_{k=l}^n\sum_{l_0+\cdots+l_k=n-k}q^{\sum_{i=0}^kil_i}\frac1{(1-q)^{n+l-k}}
s_{1,q}(k,l)(-1)^k
\frac{1}{1+q^l}.$$
where $E_{n,q}$ is the $q$-Euler numbers.

\section{$q$-Euler numbers, $q$-Stirling numbers and $q$-Bernstein polynomials
related to the fermionic $p$-adic integrals on $\mathbb Z_p$}

First, we consider the $q$-extension of the generating function of Bernstein polynomials in (\ref{AA-gen}).

For $q\in\C$ with $|1-q|_p<p^{-1/(p-1)},$ we obtain
\begin{equation}\label{B-def-ge-ft}
\begin{aligned}
F_q^{(k)}(t,x)&=\frac{t^ke^{[1-x]_qt}[x]_q^k}{k!}\\
&=[x]_q^k\sum_{n=0}^\infty \binom{n+k}{k}[1-x]_q^n\frac{t^{n+k}}{(n+k)!}\\
&=\sum_{n=k}^\infty\binom nk[x]_q^k[1-x]_q^{n-k}\frac{t^n}{n!} \\
&=\sum_{n=0}^\infty B_{k,n}(x,q)\frac{t^n}{n!},
\end{aligned}
\end{equation}
which is the generating function of the modified $q$-Bernstein type polynomials (see \cite{KCK}).
Indeed, this
generating function is also treated by Simsek and Acikgoz (see \cite{SA}).
Note that $\lim_{q\to1}F_q^{(k)}(t,x)=F^{(k)}(t,x).$
It is easy to show that
\begin{equation}\label{id-1}
[1-x]_q^{n-k}=\sum_{m=0}^\infty\sum_{l=0}^{n-k}\binom{l+m-1}{m}\binom{n-k}l(-1)^{l+m}q^l[x]_q^{l+m}(q-1)^m.
\end{equation}

From (\ref{B-op}), (\ref{de-2}), (\ref{e-p-exp}) and (\ref{id-1}), we derive the following theorem.

\begin{theorem}
For $k,n\in\mathbb Z_+$ with $n\geq k,$
$$\begin{aligned}
\begin{aligned}
\int_{\Z}\frac{B_{k,n}(x,q)}{\binom nk}&d\mu_{-1}(y) \\
&=\sum_{m=0}^\infty\sum_{l=0}^{n-k}\binom{l+m-1}{m}\binom{n-k}l(-1)^{l+m}q^l(q-1)^m
E_{l+m+k,q},
\end{aligned}
\end{aligned}$$
where $E_{n,q}$ are the $q$-Euler numbers.
\end{theorem}

It is possible to write $[x]_q^k$ as a linear combination of the modified $q$-Bernstein
polynomials by using the degree evaluation formulae and mathematical induction.
Therefore we obtain the following theorem.

\begin{theorem}[{\cite[Theorem 7]{KJY}}] \label{ber-sum}
For $k,n\in\mathbb Z_{+},i\in\mathbb N$ and $x\in[0,1],$
$$\sum_{k=i-1}^n\frac{\binom ki}{\binom ni}B_{k,n}(x,q)=[x]_q^i([x]_q+[1-x]_q)^{n-i}.$$
\end{theorem}

Let $i-1\leq n.$ Then from (\ref{Ber-def}), (\ref{id-1}) and Theorem \ref{ber-sum}, we have
\begin{equation}\label{B-def-ge}
\begin{aligned}
\,[x]_q^i&=\frac{\sum_{k=i-1}^n\frac{\binom ki\binom nk}{\binom ni}[x]_q^k[1-x]_q^{n-k}}{[x]_q^{n-i}
\left(1+\frac{[1-x]_q}{[x]_q}\right)^{n-k}} \\
&=\sum_{m=0}^\infty\sum_{k=i-1}^n\sum_{l=0}^{m+n-k}\sum_{p=0}^\infty
\frac{\binom ki\binom nk}{\binom ni}\binom{l+p-1}{p}\binom{m+n-k}{l} \\
&\quad\times\binom{n-i+m-1}{m}(-1)^{l+p+m}q^l(q-1)^p[x]_q^{i-n-m+k+p+l}.
\end{aligned}
\end{equation}
Using (\ref{q-Euler}) and (\ref{B-def-ge}), we obtain the following theorem.

\begin{theorem}\label{q-eu-Ber}
For $k,n\in\mathbb Z_{+}$ and $i\in\mathbb N$ with $i-1\leq n,$
$$\begin{aligned}
E_{i,q}
&=\sum_{m=0}^\infty\sum_{k=i-1}^n\sum_{l=0}^{m+n-k}\sum_{p=0}^\infty
\frac{\binom ki\binom nk}{\binom ni}\binom{l+p-1}{p}\binom{m+n-k}{l} \\
&\quad\times\binom{n-i+m-1}{m}(-1)^{l+p+m}q^l(q-1)^pE_{i-n-m+k+p+l,q}.
\end{aligned}$$
\end{theorem}

The $q$-String numbers of the first kind is defined by
\begin{equation}\label{1-str}
\prod_{k=1}^n(1+[k]_qz)=\sum_{k=0}^nS_1(n,k;q)z^k,
\end{equation}
and the $q$-String number of the second kind is also defined by
\begin{equation}\label{2-str}
\prod_{k=1}^n(1+[k]_qz)^{-1}=\sum_{k=0}^nS_2(n,k;q)z^k
\end{equation}
(see \cite{KCK}). Therefore, we deduce the following theorem.

\begin{theorem}[{\cite[Theorem 4]{KCK}}] \label{ber-str}
For $k,n\in\mathbb Z_{+}$ and $i\in\mathbb N,$
$$\frac{\sum_{k=i-1}^n\frac{\binom ki}{\binom ni}B_{k,n}(x,q)}{([x]_q+[1-x]_q)^{n-i}}
=\sum_{k=0}^i\sum_{l=0}^kS_1(n,l;q)S_2(i,k;q)[x]_q^l.$$
\end{theorem}

By Theorem \ref{ber-sum}, Theorem \ref{ber-str} and the definition of fermionic $p$-adic integrals on $\mathbb Z_p,$
we obtain the following theorem.

\begin{theorem} \label{Eu-str}
For $k,n\in\mathbb Z_{+}$ and $i\in\mathbb N,$
$$\begin{aligned}
E_{i,q}&=\sum_{k=i-1}^n\frac{\binom ki}{\binom ni}\int_{\Z}\frac{B_{k,n}(x,q)}{([x]_q+[1-x]_q)^{n-i}}d\mu_{-1}(x)\\
&=\sum_{k=0}^i\sum_{l=0}^kS_1(n,l;q)S_2(i,k:q)E_{l,q},
\end{aligned}$$
where $E_{i,q}$ is the $q$-Euler numbers.
\end{theorem}

Let $i-1\leq n.$ It is easy to show that
\begin{equation}\label{Ber-new}
\begin{aligned}
&[x]_q^i([x]_q+[1-x]_q)^{n-i} \\&=\sum_{l=0}^{n-i}\binom{n-i}l[x]_q^{l+i}[1-x]_q^{n-i-l} \\
&=\sum_{l=0}^{n-i}\sum_{m=0}^{n-i-l}\binom{n-i}l\binom{n-i-l}{m}(-1)^mq^m[x]_q^{m+i+l}q^{-mx} \\
&=\sum_{l=0}^{n-i}\sum_{m=0}^{n-i-l}\sum_{s=0}^\infty\binom{n-i}l\binom{n-i-l}{m}\binom{m+s-1}s
\\
&\quad\times(-1)^mq^m(1-q)^s[x]_q^{m+i+l+s}.
\end{aligned}
\end{equation}
From (\ref{Ber-new}) and Theorem \ref{ber-sum}, we have the following theorem.

\begin{theorem} \label{Eu-str-2}
For $k,n\in\mathbb Z_{+}$ and $i\in\mathbb N,$
$$\begin{aligned}
\sum_{k=i-1}^n\frac{\binom ki}{\binom ni}\int_{\Z}B_{k,n}(x,q)d\mu_{-1}(x)
&=\sum_{l=0}^{n-i}\sum_{m=0}^{n-i-l}\sum_{s=0}^\infty\binom{n-i}l\binom{n-i-l}{m}\binom{m+s-1}s
\\
&\quad\times(-1)^mq^m(1-q)^sE_{m+i+l+s,q},
\end{aligned}$$
where $E_{i,q}$ are the $q$-Euler numbers.
\end{theorem}

In the same manner, we can obtain the following theorem.

\begin{theorem} \label{Eu-str-3}
For $k,n\in\mathbb Z_{+}$ and $i\in\mathbb N,$
$$\begin{aligned}
\int_{\Z}B_{k,n}(x,q)d\mu_{-1}(x)
&=\sum_{j=k}^n\sum_{m=0}^\infty\binom jk\binom nj\binom{j-k+m-1}{m}
\\
&\quad\times(-1)^{j-k+m}q^{j-k}(q-1)^mE_{m+j,q},
\end{aligned}$$
where $E_{i,q}$ are the $q$-Euler numbers.
\end{theorem}

\section{Further remarks and observations}

The $q$-binomial formulas are known,
\begin{equation}\label{q-bi-f-two}
\begin{aligned}
&(a;q)_n=(1-a)(1-aq)\cdots(1-aq^{n-1})=\sum_{i=0}^n\binom ni_q q^{\binom i2}(-1)^ia^i, \\
&\frac1{(a;q)_n}=\frac1{(1-a)(1-aq)\cdots(1-aq^{n-1})}=\sum_{i=0}^\infty\binom{n+i-1}{i}_qa^i.
\end{aligned}
\end{equation}
For $h\in\mathbb Z,n\in\mathbb Z_+$ and $r\in\mathbb N,$
we introduce the extended higher-order $q$-Euler polynomials as follows \cite{KT09}:
\begin{equation}\label{h-o-q-Eu}
E_{n,q}^{(h,r)}(x)=\int_{\Z}\cdots\int_{\Z}
q^{\sum_{j=1}^r(h-j)x_j}[x+x_1+\cdots+x_r]_q^nd\mu_{-1}(x_1)\cdots d\mu_{-1}(x_r).
\end{equation}
Then
\begin{equation}\label{h-o-ex}
\begin{aligned}
E_{n,q}^{(h,r)}(x)&=\frac{2^r}{(1-q)^n}\sum_{l=0}^n\binom nl (-1)^l\frac{q^{lx}}{(-q^{h-1+l};q^{-1})_r} \\
&=\frac{2^r}{(1-q)^n}\sum_{l=0}^n\binom nl (-1)^l\frac{q^{lx}}{(-q^{h-r+l};q)_r}.
\end{aligned}
\end{equation}
Let us now define the extended higher-order N\"orlund type $q$-Euler polynomials as follows \cite{KT09}:
\begin{equation}\label{h-o-ex-2}
\begin{aligned}
E_{n,q}^{(h,-r)}(x)&=\frac1{(1-q)^n}\sum_{l=0}^n\binom nl(-1)^l \\
&\quad\times\frac{q^{lx}}{\int_{\Z}\cdots\int_{\Z}
q^{l(x_1+\cdots+x_r)}q^{\sum_{j=1}^r(h-j)x_j}d\mu_{-1}(x_1)\cdots d\mu_{-1}(x_r)}.
\end{aligned}
\end{equation}
In the special case $x=0,$ $E_{n,q}^{(h,-r)}=E_{n,q}^{(h,-r)}(0)$ are called
the extended higher-order N\"orlund type $q$-Euler numbers.
From (\ref{h-o-ex-2}), we note that
\begin{equation}\label{h-o-ex-3}
\begin{aligned}
E_{n,q}^{(h,-r)}(x)&=\frac{1}{2^r(1-q)^n}\sum_{l=0}^n\binom nl (-1)^lq^{lx}(-q^{h-r+l};q)_r\\
&=\frac1{2^r}\sum_{m=0}^rq^{\binom m2}q^{(h-r)m}\binom rm_q[m+x]_q^n.
\end{aligned}
\end{equation}
A simple manipulation shows that
\begin{equation}\label{st-1}
q^{\binom m2}\binom rm_q =\frac{q^{\binom m2}[r]_q\cdots[r-m+1]_q}{[m]_q!}=\frac1{[m]_q!}
\prod_{k=0}^{m-1}([r]_q-[k]_q)
\end{equation}
and
\begin{equation}\label{st-2}
\prod_{k=0}^{n-1}(z-[k]_q)=z^n\prod_{k=0}^{n-1}\left(1-\frac{[k]_q}{z}\right)
=\sum_{k=0}^n S_1(n-1,k;q)(-1)^kz^{n-k}.
\end{equation}
Formulas (\ref{h-o-ex-3}), (\ref{st-1}) and (\ref{st-2}) imply the following lemma.

\begin{lemma} \label{Eu-str-pro}
For $h\in\mathbb Z,n\in\mathbb Z_+$ and $r\in\mathbb N,$
$$E_{n,q}^{(h,-r)}(x)=\frac{1}{2^r[m]_q!}\sum_{m=0}^r\sum_{k=0}^mq^{(h-r)m}S_1(m-1,k;q)(-1)^k[r]_q^{m-k}
[x+m]_q^n.$$
\end{lemma}

From (\ref{q-id-1}), we can easily see that
\begin{equation}\label{le-ber-id}
[x+m]_q^n=\frac1{(1-q)^n}\sum_{j=0}^n\sum_{l=0}^j\binom nj\binom jl(-1)^{j+l}(1-q)^lq^{mj}[x]_q^l.
\end{equation}
Using (\ref{q-Euler}) and (\ref{le-ber-id}), we obtain the following lemma.

\begin{lemma} \label{Eu-str-pro2}
For $m,n\in \mathbb Z_+,$
$$E_{n,q}(m)=\frac1{(1-q)^n}\sum_{j=0}^n\sum_{l=0}^j\binom nj\binom jl(-1)^{j+l}(1-q)^lq^{mj}E_{l,q}.$$
\end{lemma}

By Lemma \ref{le-ber-id}, lemma \ref{Eu-str-pro2} and the definition of fermionic $p$-adic integrals on $\mathbb Z_p,$
we obtain the following theorem.

\begin{theorem} \label{mul-q-Eu}
For $h\in\mathbb Z,n\in\mathbb Z_+$ and $r\in\mathbb N,$
$$\begin{aligned}
\int_{\Z}E_{n,q}^{(h,-r)}(x)d\mu_{-1}(x)&=\frac{2^{-r}}{[m]_q!}\sum_{m=0}^r\sum_{k=0}^mq^{(h-r)m}S_1(m-1,k;q)(-1)^k[r]_q^{m-k}
E_{n,q}(m) \\
&=\frac{1}{2^r[m]_q!}\sum_{m=0}^r\sum_{k=0}^mq^{(h-r)m}S_1(m-1,k;q)(-1)^k[r]_q^{m-k} \\
&\quad\times\frac1{(1-q)^n}\sum_{j=0}^n\sum_{l=0}^j\binom nj\binom jl(-1)^{j+l}(1-q)^lq^{mj}E_{l,q}.
\end{aligned}
$$
\end{theorem}

Put $h=0$ in (\ref{h-o-ex-2}). We consider the following polynomials $E_{n,q}^{(0,-r)}(x):$
\begin{equation}\label{h-o-ex-h=0}
\begin{aligned}
E_{n,q}^{(0,-r)}(x)&=\sum_{l=0}^n\frac{(1-q)^{-n}\binom nl(-1)^lq^{lx}}{\int_{\Z}\cdots\int_{\Z}
q^{l(x_1+\cdots+x_r)}q^{-\sum_{j=1}^rjx_j}d\mu_{-1}(x_1)\cdots d\mu_{-1}(x_r)}.
\end{aligned}
\end{equation}
Then
$$E_{n,q}^{(0,-r)}(x)=\frac1{2^r}\sum_{m=0}^r\binom rmq^{\binom m2-rm}[m+x]_q^n.$$
A simple calculation of the ferminionic $p$-adic invariant integral on $\Z$ shows that
$$\int_{\Z}E_{n,q}^{(0,-r)}(x)d\mu_{-1}(x)=\frac1{2^r}\sum_{m=0}^r\binom rmq^{\binom m2-rm}E_{n,q}(m).$$
Using Theorem \ref{mul-q-Eu}, we can also prove that
$$\int_{\Z}E_{n,q}^{(0,-r)}(x)d\mu_{-1}(x)=\frac{2^{-r}}{[m]_q!}\sum_{m=0}^r\sum_{k=0}^mq^{-rm}S_1(m-1,k;q)(-1)^k[r]_q^{m-k}
E_{n,q}(m).$$
Therefore, we obtain the following theorem.

\begin{theorem}
For $m\in\mathbb Z_+,r\in\mathbb N$ with $m\leq r,$
$$\binom rmq^{\binom m2-rm}
=\frac{1}{[m]_q!}\sum_{k=0}^mq^{-rm}S_1(m-1,k;q)(-1)^k[r]_q^{m-k}.$$
\end{theorem}

\bibliography{central}

\begin{thebibliography}{00}

\bibitem{AA} M. Acikgoz and S. Araci,
\textit{A study on the integral of the product of several type Bernstein polynomials},
IST Transaction of Applied Mathematics-Modelling and Simulation, 2010.

\bibitem{AA2} M. Acikgoz and S. Araci,
\textit{On the generating function of the Bernstein polynomials},
Accepted to AIP on 24 March 2010 for ICNAAM 2010.

\bibitem{Be} S. Bernstein,
\textit{Demonstration du theoreme de Weierstrass, fondee sur le calcul des probabilities},
Commun. Soc. Math. Kharkow (2) \textbf{13} (1912-1913), 1--2.

\bibitem{GG} N. K. Govil and V. Gupta,
\textit{Convergence of $q$-Meyer-Konig-Zeller-Purrmeyer operators},
Adv. Stud. Contemp. Math. \textbf{19} (2009), 181--189.

\bibitem{KT1} T. Kim,
\textit{On a $q$-analogue of the $p$-adic log gamma functions and related integrals},
J. Number Theory \textbf{76} (1999), 320--329.

\bibitem{KT2} T. Kim,
\textit{$q$-Volkenborn integration},
Russ. J. Math. Phys. \textbf{9} (2002), no. 3, 288--299.

\bibitem{KT3} T. Kim,
\textit{Power series and asymptotic series associated with the $q$-analog of the two-variable $p$-adic $L$-function},
Russ. J. Math. Phys. \textbf{12} (2005), 186--196.

\bibitem{KT4} T. Kim,
\textit{On the analogs of Euler numbers and polynomials associated with $p$-adic $q$-integral on $\mathbb Z_p$ at $q=-1$},
J. Math. Anal. Appl. \textbf{331} (2007), 779--792.

\bibitem{KT5} T. Kim,
\textit{$q$-Bernoulli numbers and polynomials associated with Gaussian binomial coefficients},
Russ. J. Math. Phys. \textbf{15} (2008), 51--57.

\bibitem{KT6} T. Kim,
\textit{On $p$-adic interpolating function for $q$-Euler numbers and its derivatives},
J. Math. Anal. Appl. \textbf{339} (2008), 598--608.

\bibitem{KT09} T. Kim,
\textit{Some identities on the $q$-Euler polynomials of higher order and $q$-Stirling numbers
by the fermionic $p$-adic integral on $\mathbb Z_p$},
Russ. J. Math. Phys. \textbf{16} (2009), no. 4, 484--491.

\bibitem{KT7} T. Kim,
\textit{$q$-Bernoulli numbers associated with $q$-Stirling numbers},
Adv. Difference Equ. 2008, Art. ID 743295, 10 pp.

\bibitem{KT8} T. Kim,
\textit{Barnes-type multiple $q$-zeta functions and $q$-Euler polynomials},
J. Physics A : Math. Theor., \textbf{43} (2010), 255201.

\bibitem{KJY} T. Kim, L.-C. Jang and H. Yi,
\textit{Note on the modified $q$-Bernstein polynomials}, Discrete Dynamics in Nature and Society (in Press).
arXiv 1005.4293.

\bibitem{KCK} T. Kim, J. Choi and Y.-H. Kim
\textit{Some identities on the $q$-Bernstein polynomials, $q$-Stirling numbers and $q$-Bernoulli numbers},
Adv. Stud. Contemp. Math. 20 (2010), no. 3, 335--341.
arXiv 1006.2033.

\bibitem{Ph} G. M. Phillips,
\textit{Bernstein polynomials based on the $q$-integer},
Annals of Numerical Analysis \textbf{4} (1997), 511--518.

\bibitem{SA}
Y. Simsek and M. Acikgoz, \textit{A new generating function of q-Bernstein-type polynomials and their interpolation function},
 Abstr. Appl. Anal. ID 769095 (2010), 12 pp.



\end{thebibliography}

\end{document}